\documentclass[twocolumn,10pt]{article} 
\usepackage{times}
\usepackage{amsmath,amssymb,amscd}
\usepackage{epic}
\usepackage{ecltree}
\usepackage{theorem}
\usepackage{ifthen}

\theoremheaderfont{\scshape}
\theorembodyfont{\normalfont\slshape}
\theoremstyle{plain}
\newtheorem{theorem}{Theorem}
\newtheorem{proposition}{Proposition}
\newtheorem{lemma}{Lemma}

\newtheorem{definition}{Definition}

\theorembodyfont{\normalfont}

\makeatletter
\renewcommand{\section}{\@startsection{section}{1}{0mm}%
                                   {- \baselineskip}%
                                   {.5\baselineskip}%
                                   {\normalfont\large\bf}}
\renewcommand\subsection{\@startsection{subsection}{2}{\z@}%
                                     {-3.25ex\@plus -1ex \@minus -.2ex}%
                                     {1.5ex \@plus .2ex}%
                                     {\normalfont\normalsize\bfseries}}
\renewcommand\subsubsection{\@startsection{subsubsection}{3}{\z@}%
                                     {-3.25ex\@plus -1ex \@minus -.2ex}%
                                     {1.5ex \@plus .2ex}%
                                     {\normalfont\normalsize\bfseries}}
\renewcommand\paragraph{\@startsection{paragraph}{4}{\z@}%
                                    {3.25ex \@plus1ex \@minus.2ex}%
                                    {-1em}%
                                    {\normalfont\normalsize\bfseries}}
\renewcommand\subparagraph{\@startsection{subparagraph}{5}{\parindent}%
                                       {3.25ex \@plus1ex \@minus .2ex}%
                                       {-1em}%
                                      {\normalfont\normalsize\bfseries}}
\makeatother

\newenvironment{proof}
  {\begin{trivlist}\item{}\normalfont\textit{Proof. }}
  {\hfill$\square$\end{trivlist}}
\newenvironment{proofof}[1]
  {\begin{trivlist}\item{}\normalfont\textit{Proof\, (of #1). }}
  {\hfill$\square$\end{trivlist}}
\newcommand{\defn}[1]{{\textit{\textbf{#1}}}}
\newcommand{\ie}{\textit{i.e.}}
\newcommand{\Ie}{\textit{I.e.}}
\newcommand{\eg}{\textit{e.g.}}

\newcommand{\etc}{\textit{etc.}}

\newlength{\textw}
\newcounter{boxlabelcounter}

\newcommand{\boxref}[1]{\the#1}

\def\POl{}%
\newcount\PLv\newcount\PLw\newcount\PLx\newcount\PLy\newdimen\PLyy\newdimen\PLX
\newbox\PLdot\setbox\PLdot\hbox{{\tiny.}}\def\scl{.07}% resettable scale
\def\PLot#1{\PLx`#1\advance\PLx-42\PLy\PLx\PLv\PLx\divide\PLy9\PLw\PLy\multiply
\PLw9\advance\PLx-\PLw\advance\PLx-4\PLy-\PLy\advance\PLy4\PLX=\POl\the\PLx pt
\advance\PLyy\POl\the\PLy pt\wd\PLdot=\scl\PLX\raise\scl\PLyy\copy\PLdot}
\def\draw#1{\ifx#1\end\let\next=\relax\else\PLot#1\let\next=\draw\fi\next}

\newcommand{\longto}{\longrightarrow}

\renewcommand{\perp}{^\bot}

\newcommand{\falseconst}{\mathsf{false}}
\newcommand{\trueconst}{\mathsf{true}}

\usepackage{calc}
\newcounter{l}
\newcommand{\rln}[3]{\setcounter{l}{#2-#1}\put(#1,#3){\line(1,0){\value{l}}}}% x1 x2 y   -
\newcommand{\lln}[3]{\setcounter{l}{#1-#2}\put(#2,#3){\line(1,0){\value{l}}}}% x1 x2 y   -
\newcommand{\uln}[3]{\setcounter{l}{#3-#2}\put(#1,#2){\line(0,1){\value{l}}}}% x y1 y2   |
\newcommand{\dln}[3]{\setcounter{l}{#2-#3}\put(#1,#2){\line(0,-1){\value{l}}}}% x y1 y2   |
\newcommand{\hln}[3]{\ifthenelse{#2 > #1}%  x1 x2 y  -
{\rln{#1}{#2}{#3}}{\lln{#1}{#2}{#3}}}
\newcommand{\vln}[3]{\ifthenelse{#3 > #2}%  x y1 y2   |
{\uln{#1}{#2}{#3}}{\dln{#1}{#2}{#3}}}

\newlength{\overlaylen}

\newcommand{\id}{\mathsf{id}}

\newcommand{\myitem}[1]{\item[(#1)]}

\newcommand{\edgelen}{30}
\newcommand{\edgerad}{15}
\newcommand{\putmathlabel}[3]{\put(#1,#2){\makebox(0,0){$#3$}}}

\newcommand{\crossgraph}{%
\put(-\edgerad,-\edgerad){%
\putmathlabel{0}{\edgelen}{p\;\;}%
\putmathlabel{0}{0}{q\;\;}%
\putmathlabel{\edgelen}{\edgelen}{\;\;p}%
\putmathlabel{\edgelen}{0}{\;\;\overline{p}}%
\put(3,30){\line(1,0){24}}%
\put(3,27){\line(1,-1){24}}%
\put(3,0){\line(1,0){24}}%
\put(3,3){\line(1,1){24}}}}
\newcommand{\crossresolutions}{%
\put(-\edgerad,-\edgerad){%
\put(0,30){\oval(20,20)[t]}
\put(30,30){\oval(20,20)[t]}
\put(-10,0){\line(0,1){30}}
\put(10,0){\line(0,1){30}}
\put(20,0){\line(0,1){30}}
\put(40,0){\line(0,1){30}}
\put(0,0){\oval(20,20)[b]}
\put(30,0){\oval(20,20)[b]}
\putmathlabel{0}{\edgelen}{p}%
\putmathlabel{0}{0}{q}%
\putmathlabel{\edgelen}{\edgelen}{p}%
\putmathlabel{\edgelen}{0}{\overline{p}}}}
\newcommand{\mapstowidth}{65}

\newcommand{\G}{\mathalpha\mathbf{G\hspace{.05ex}}}
\newcommand{\Rel}{\mathsf{Rel}}
\newcommand{\grel}[1]{\G(\Rel^{#1})}

\newcommand{\laxg}{\G_{{}^\le}}
\newcommand{\glrel}[1]{\laxg(\Rel^{#1})}
\newcommand{\lgrel}[1]{\glrel{#1}}

\newcommand{\inj}{\iota}

\newcommand{\cwk}{\mathsf{CW}_K{}}

\newcommand{\dlk}{\mathsf{D}_K}
\newcommand{\res}{($\mathcal{R}$)}
\newcommand{\resprime}{($\mathcal{R}'$)}
\newcommand{\lres}{($\mathcal{R}_{{}^\le}$)}
\newcommand{\lresprime}{($\mathcal{R}'_{{}^\le}$)}
\newcommand{\sqleft}[1]{\put(0,-#1){\line(-1,0){#1}}\put(-#1,0){\line(0,1){#1}\line(0,-1){#1}}\put(0,#1){\line(-1,0){#1}}}
\newcommand{\sqright}[1]{\put(0,#1){\line(1,0){#1}\line(0,-1){#1}}\put(0,-#1){\line(1,0){#1}\line(0,1){#1}}}
\newcommand{\uax}{\mbox{\footnotesize{$\mathsf{AX}$}}}
\newcommand{\ucut}{\mbox{\footnotesize{$\mathsf{CUT}$}}}

\newcommand{\bul}{\mathsf{B}^u_L}
\newcommand{\bll}{\mathsf{B}^a_L}
\newcommand{\blinkingl}{\mathsf{B}^l_L}
\newcommand{\Link}{\mathsf{Link}}
\newcommand{\dual}[1]{\overline{#1}}
\newcommand{\putlabel}[3]{\put(#1,#2){\makebox(0,0){#3}}}

\newcommand{\bulle}{\raisebox{0pt}[5.8pt][0pt]{\large$\bullet$}}
\newcommand{\putv}[2]{\putlabel{#1}{#2}{\bulle}}
\newcommand{\assoc}{\mathsf{assoc}}
\title{
\small
\vspace{-7ex}
{\bf\Large Logic Without Syntax
\\\vspace{-2ex}%
}%
\author{\normalsize
   \sc Dominic Hughes \\ \normalsize Stanford University \\[-.3ex]
   \footnotesize
\vspace{-3ex}
}}%
\date{}
\begin{document}

\maketitle
%\begin{abstract}
%\paragraph*{Abstract}
\sl This paper presents an abstract, mathematical formulation of classical
propositional logic.
It proceeds layer by layer: (1) abstract, syntax-free propositions;
(2) abstract, syntax-free contraction-weakening proofs;
(3) distribution; (4) axioms $p\vee\overline{p}$.

Abstract propositions correspond to objects of the category
$\G(\Rel^{L})$ where $\G$ is the Hyland-Tan double glueing
construction, $\Rel$ is the standard category of sets and relations,
and $L$ is a set of literals.
%and cuts $\overline{p}\wedge p$.
%\vspace{-1.5ex}%
%Blah blah blah
%
Abstract proofs are morphisms of a tight orthogonality subcategory of
$\lgrel{L}$, where we define $\laxg$ as a lax variant of $\G$.
We prove that the free binary product-sum category
(contraction-weakening logic) over $L$ is a full subcategory of
$\grel{L}$, and the free distributive lattice category
(contraction-weakening-distribution logic) is a full subcategory of
$\lgrel{L}$.
We explore general constructions for adding axioms, which are not
$\Rel$-specific or ($p\vee\dual{p}$)-specific.
%\end{abstract}
%\vspace{-4ex}
\normalfont
\section{Introduction}\label{intro}
%\vspace{-2.5ex}
\paragraph*{Abstract propositions.} Typically logicians define a 
proposition or formula as a labelled tree.  Using de Morgan duality
($\neg(A\wedge B)=\neg A\vee \neg B$ \etc) one needs only trees
labelled by literals (variables and their duals) and constants on
leaves and $\vee$ and $\wedge$ on internal nodes.  To quotient by
associativity and commutivity, one may use graphs (cographs or
$P_4$-free graphs \cite{CLS81}), for example,
\begin{center}
\begin{picture}(0,22)
\thicklines\put(0,12){\putmathlabel{-\mapstowidth}{0}{(p\vee q)\wedge(p\vee \overline{p})}
\putmathlabel{0}{0}{\mapsto}
\put(\mapstowidth,0){\crossgraph}}
\end{picture}
\end{center}
drawing an edge between leaves iff they meet in the parse tree at a
$\wedge$.  In this paper we go a step further, and define an
\emph{abstract proposition} as a set of \emph{leaves}
together with a set of subsets, called \emph{resolutions}.  For
example,
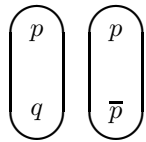
\begin{center}
\begin{picture}(0,42)
\thicklines\put(0,21){\crossresolutions}
\end{picture}
\end{center}
has the four leaves of the formula/graph depicted earlier, and two
resolutions, the maximal independent sets (maximal co-cliques) of the
graph.  
%Conjunction $\wedge$ is disjoint union.  
Any syntactic proposition can be reconstructed from its
abstract leaf/resolution presentation.
(The terminology `resolution' 
%Resolutions (in syntactic guise)
%, and the \emph{resolution condition},
here comes from the definition of MALL proof net \cite{HG03}.)
%[The term `resolution'  from the \emph{resolution condition}
%for MALL proof nets \cite{??}.]

A key advantage of this abstraction is a crisp mathematical treatment
of the logical units/constants false $0$ and true $1$.  In the traditional
syntactic world, $0$ and $1$ have the same stature as literals, taking
up actual ink on the page as labelled leaves.  They are artificially
dual to each other, by \emph{fiat}, and artificially act as units for
syntactic $\vee$ and $\wedge$.  The graphical representation gets
closer to a nice treatment of units: the empty graph $\epsilon$ is a
unit for the operations $\vee$ and $\wedge$ (union and join) on
graphs; however then one has degeneracy $\epsilon=0=1$, so one must
resort once again to artificially promoting $0$ and $1$ to actual
labelled vertices, and the units remain ad hoc.

Abstract propositions succeed in having both units empty (no leaves),
hence mathematically crisp as units for the operations $\vee$ and
$\wedge$, \emph{without} identifying them (as in the graph case
$\epsilon=0=1$):
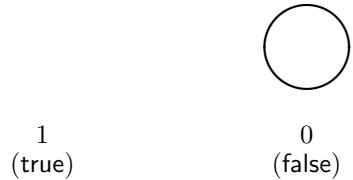
\begin{center}
\begin{picture}(0,65)(0,-18)
\thicklines\put(50,30){\circle{30}}\putmathlabel{-50}{-3}{1}
\putmathlabel{-50}{-15}{(\trueconst)}
\putmathlabel{50}{-3}{0}
\putmathlabel{50}{-15}{(\falseconst)}
\end{picture}
\end{center}
The unit $1$ has no leaves and no resolution, and $0$ has no leaves
with one (empty) resolution.

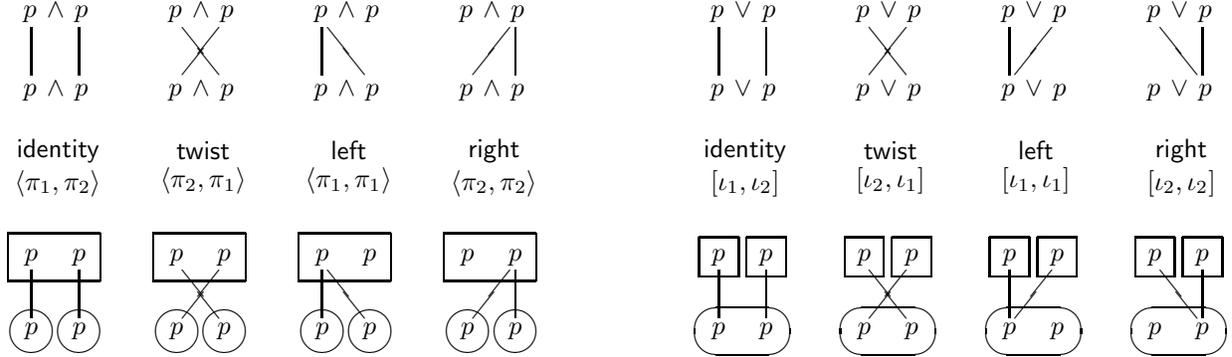
\begin{figure*}[t]
\newcommand{\aatop}{26}\newcommand{\aabot}{8}
\newcommand{\underlabel}[2]{\put(13,-27){\makebox(0,0){\shortstack{$#1$\\$#2$}}}}
\newcommand{\ps}{\put(0,0){$p\,\wedge\, p$}\put(0,30){$p\,\wedge\, p$}}
\newcommand{\plines}[4]{\drawline(#1,\aatop)(#2,\aabot)\drawline(#3,\aatop)(#4,\aabot)}
\newcommand{\aapic}[4]{\ps\plines{#1}{#2}{#3}{#4}}
% FOR OLD FORM, SHOWING BOTH CLIQUES AND CO-CLIQUES
%\newcommand{\diam}{14}\newcommand{\rad}{7}
%\newcommand{\bigdiam}{22}\newcommand{\bigrad}{11}
%\newcommand{\sqdiam}{14}\newcommand{\sqrad}{7}
%\newcommand{\bigsqdiam}{20}\newcommand{\bigsqrad}{10}
%\newcommand{\width}{18}
%\newcommand{\bigsqwidth}{38}
%\newcommand{\poval}{\makebox(0,0){$p$}\oval(\diam,\diam)}
%\newcommand{\psquare}{\makebox(0,0){$p$}%
%\put(-\sqrad,-\sqrad){\line(1,0){\sqdiam}\line(0,1){\sqdiam}}\put(\sqrad,\sqrad){\line(-1,0){\sqdiam}\line(0,-1){\sqdiam}}}
%\newcommand{\ppoval}{%
%%\makebox(0,0){$p$}\put(18,0){\makebox(0,0){$p$}}
%\put(0,0){\psquare}\put(\width,0){\psquare}%
%\oval(\bigdiam,\bigdiam)[l]\put(0,\bigrad){\line(1,0){\width}}\put(0,-\bigrad){\line(1,0){\width}}\put(\width,0){\oval(\bigdiam,\bigdiam)[r]}}
%\newcommand{\ppsquare}{%
%\put(0,0){\poval}\put(\width,0){\poval}%
%\put(-\bigsqrad,-\bigsqrad){\line(1,0){\bigsqwidth}\line(0,1){\bigsqdiam}}%
%\put(-\bigsqrad,\bigsqrad){\line(0,-1){\bigsqdiam}\line(1,0){\bigsqwidth}}}
% END OF OLD FORM
\newcommand{\diam}{16}\newcommand{\rad}{8}
\newcommand{\bigdiam}{18}\newcommand{\bigrad}{9}
\newcommand{\sqdiam}{15}\newcommand{\sqrad}{7.5}
\newcommand{\bigsqdiam}{18}\newcommand{\bigsqrad}{9}
\newcommand{\width}{18}
\newcommand{\bigsqwidth}{35}
\newcommand{\poval}{\makebox(0,0){$p$}\oval(\diam,\diam)}
\newcommand{\psquare}{\makebox(0,0){$p$}%
\put(-\sqrad,-\sqrad){\line(1,0){\sqdiam}\line(0,1){\sqdiam}}\put(\sqrad,\sqrad){\line(-1,0){\sqdiam}\line(0,-1){\sqdiam}}}
\newcommand{\ppoval}{%
\put(0,0){\makebox(0,0){$p$}}\put(\width,0){\makebox(0,0){$p$}}%
\oval(\bigdiam,\bigdiam)[l]\put(0,\bigrad){\line(1,0){\width}}\put(0,-\bigrad){\line(1,0){\width}}\put(\width,0){\oval(\bigdiam,\bigdiam)[r]}}
\newcommand{\ppsquare}{%
\put(0,0){\makebox(0,0){$p$}}\put(\width,0){\makebox(0,0){$p$}}%
\put(-\bigsqrad,-\bigsqrad){\line(1,0){\bigsqwidth}\line(0,1){\bigsqdiam}}%
\put(-\bigsqrad,\bigsqrad){\line(0,-1){\bigsqdiam}\line(1,0){\bigsqwidth}}}
\newcommand{\psquares}{\put(0,0){\psquare}\put(\width,0){\psquare}}
\newcommand{\povals}{\put(0,0){\poval}\put(\width,0){\poval}}
\begin{center}%
\begin{picture}(0,139)
\put(-220,100){%
\put(0,0){\underlabel{\mathsf{identity}}{\langle \pi_1,\pi_2\rangle}\aapic{3}{3}{21}{21}}
\put(55,0){\underlabel{\mathsf{twist}}{\langle \pi_2,\pi_1\rangle}\aapic{5}{19}{19}{5}}
\put(110,0){\underlabel{\mathsf{left}}{\langle \pi_1,\pi_1\rangle}\aapic{3}{3}{5}{19}}
\put(165,0){\underlabel{\mathsf{right}}{\langle \pi_2,\pi_2\rangle}\aapic{19}{5}{21}{21}}}
\put(0,-100){}
\renewcommand{\ps}{\put(3,3){\povals}\put(3,31){\ppsquare}}
\put(-220,8){%
\put(0,0){\aapic{3}{3}{21}{21}}
\put(55,0){\aapic{5}{19}{19}{5}}
\put(110,0){\aapic{3}{3}{5}{19}}
\put(165,0){\aapic{19}{5}{21}{21}}}
\renewcommand{\ps}{\put(0,0){$p\,\vee\, p$}\put(0,30){$p\,\vee \,p$}}
\put(40,100){%
\put(0,0){\underlabel{\mathsf{identity}}{[ \inj_1,\inj_2]}\aapic{3}{3}{21}{21}}
\put(55,0){\underlabel{\mathsf{twist}}{[ \inj_2,\inj_1]}\aapic{5}{19}{19}{5}}
\put(110,0){\underlabel{\mathsf{left}}{[ \inj_1,\inj_1]}\aapic{3}{3}{19}{5}}
\put(165,0){\underlabel{\mathsf{right}}{[ \inj_2,\inj_2]}\aapic{5}{19}{21}{21}}}
\renewcommand{\ps}{\put(3,3){\ppoval}\put(3,31){\psquares}}
\put(40,8){%
\put(0,0){\aapic{3}{3}{21}{21}}
\put(55,0){\aapic{5}{19}{19}{5}}
\put(110,0){\aapic{3}{3}{19}{5}}
\put(165,0){\aapic{5}{19}{21}{21}}}
\end{picture}
\end{center}
\caption{\label{projinj-fig}The four abstract proofs $p\wedge p\to p\wedge p$, and dually,  $p\vee p\to p\vee p$.
Each proof is a $\G(\Rel^L)$ morphism, a binary relation between
leaves which satisfies the resolution condition.
The top row shows the eight morphisms in syntactic form. The bottom
row shows the same morphisms between the corresponding abstract
propositions, where the target propositions are specified by their
resolutions (curved regions), and the source propositions by their
coresolutions (rectangular regions).
Note that the resolution condition is satisfied: 
%(a) resolutions
%(curved regions) map upwards to resolutions and coresolutions
%(rectangular) map downwards to coresolutions, or equivalently, (b) 
there is a unique edge between any output resolution (curved region) and input coresolution
(square region).
The resolution condition characterises free binary product-sum
categories.%
%Each morphism is a canonical map in the free product-sum category.
%
}\end{figure*}

\paragraph*{Abstract proofs.}  
%We define abstract proofs as morphisms between abstract propositions.
%
Abstract propositions correspond to certain objects of the category
$\G(\Rel^{L})$ studied in \cite{Hug04}, where
\begin{itemize} 
\item $\G$ is the Hyland-Tan double glueing
construction \cite{Tan97},
\item $\Rel$ is the standard category of sets and
relations, and
\item $L=\{p,\overline{p},q,\overline{q},\ldots\}$ is a set of literals.
%$2V=V+\overline{V}$, for $V$ a set of propositional
%variables.
\end{itemize}
Thus a $\G(\Rel^{L})$ morphism $A\to B$ provides an off-the-shelf
notion of an abstract proof of $B$ from $A$.  By definition, a
$\grel{L}$ morphism is a certain kind of binary relation between the
leaves of $A$ and the leaves of $B$.
Figure~\ref{projinj-fig} shows the four morphisms
$p\wedge p\to p\wedge p$, and dually, the four morphisms $p\vee p\to
p\vee p$.
%which respects labels and satisfies either of the following equivalent
%conditions:
By definition of double glueing, a morphism $R$ must satisfy:
\begin{itemize}
\item[\res] \emph{Resolution condition.} $R$ pulls resolutions backwards and
pushes coresolutions forwards.
%\footnote{For logical reaons, we
%have taken resolutions as the contravariant part, so the ``co'' is
%opposite to usual $\G(\mathcal{C})$.}
\end{itemize}
More precisely, $R:A\to B$ maps resolutions of $B$ to resolutions of
$A$, and coresolutions of $A$ to coresolutions of $B$, where a
coresolution of $X$ is a resolution of its dual $\neg X$.
In the special case that $A$ and $B$ are syntactic
propositions, this coincides with the usual resolution condition on
MALL proof nets \cite{HG03}:
%\begin{itemize}
%\item \emph{Resolution condition.} $R$ has a unique edge on any resolution of $A\to B$.
%\end{itemize}
%Here a resolution of $A\to B$ is a resolution of $\overline{A}\vee B$,
%\ie, a union of a coresolution of $A$ and a resolution of $B$.
%In other words:
\begin{itemize}
\item[\resprime] \emph{Resolution condition.} $R$ has a unique edge between any 
output resolution and input coresolution.\footnote{\label{mall}Recall from
\cite{HG03} the resolution condition on a set $R$ of linkings
on a sequent or proposition $\Gamma$: \emph{$R$ has a unique linking
on any resolution of $\Gamma$.}
%$R$ has a unique edge on any resolution of $A\to B$.
In the current pure additive (\ie\ product-sum) 
%$\G(\Rel^L)$ 
setting every linking has just one edge, \ie, $R$ is simply a set of
edges.  The main text quotes this condition with $\Gamma\,=\,A\to
B\,=\,\overline{A}\vee B$. A resolution of $\overline{A}\vee B$ is a
union of a coresulution of $A$ and a resolution of $B$.}
\end{itemize}
%It may be instructive to verify this condition for the examples in
%Figure~\ref{projinj-fig}.  
We prove:
%\footnote{In a more syntactic form, I noticed that the MALL
%resolution condition characterises the free product-sum category a few years ago.
\begin{itemize}
\item[{}] \emph{Theorem.} 
$\G(\Rel^K)$ contains as a full subcategory the free binary product-sum
category generated by the set $K$.
%The full subcategory $\cwk$ of $\G(\Rel^K)$ whose objects are
%syntactic propositions is the free product-sum category generated by
%the set $K$.
%free product-sum category generated by $V$ is a full subcategory of
%$\lzero$.
% $\Lprodsum$.
\end{itemize}
%where an abstract proposition is \emph{syntactic} if it is finite and
%non-empty.  
Thus 
%$\cwk$ is 
we obtain an abstract, syntax-free formulation of pure
contraction-weakening logic over a set of atoms $K$: every morphism
(abstract proof) $A\to B$ is a composite of the units of the
product/sum adjunctions, the natural transformations (inferences)
%projection, injection, diagonal and codiagonal maps,
\begin{displaymath}
\begin{array}{cccccl}
\pi_i &:& A_1\wedge A_2 & \to & A_i & (\text{projection})\\
\inj_i &:& A_i & \to & A_1\vee A_2 & (\text{injection})\\
\delta &:& A & \to & A\wedge A & (\text{diagonal})\\
\epsilon &:& A\vee A & \to & A & (\text{codiagonal})
%\begin{array}{ccccc@{\hspace{7ex}}ccccc}
%\inj_i &:& A_i & \to & A_1\vee A_2 & \delta &:& A & \to & A\wedge A\\
%\pi_i &:& A_1\wedge A_2 & \to & A_i & \epsilon &:& A\vee A & \to & A.
\end{array}
\end{displaymath}
%with equality between proofs given by the equalities that hold in every product-sum ca
\paragraph*{Adding distribution.} The obvious candidates for a distribution 
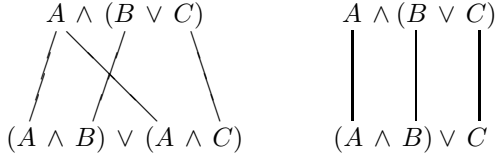
\begin{center}
\begin{picture}(0,60)
\put(-55,0){%
\putmathlabel{0}{50}{A\,\wedge\,(B\,\vee\,C)}
\drawline(-26,44)(-36,10)
\drawline(-22,44)(12,10)
\drawline(0,44)(-12,10)
\drawline(25,44)(36,10)
\putmathlabel{0}{3}{(A\,\wedge\,B)\,\vee\,(A\,\wedge\,C)}}
\put(55,0){%
\putmathlabel{0}{50}{\;A\,\wedge(B\,\vee\,C)}
\drawline(-24,44)(-24,10)
\drawline(0,44)(0,10)
\drawline(24,44)(24,10)
\putmathlabel{-2}{3}{(A\,\wedge\,B)\vee\,C}}
\end{picture}
\end{center}
fail the resolution condition.  They are not $\G(\Rel^L)$ morphisms.
Condition \res\ fails because the image of an output
resolution is strictly larger than an input resolution:
%
%\begin{center}
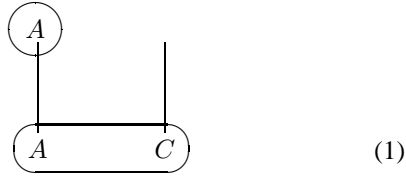
\begin{equation}\label{firstfail}
\newcommand{\sqdiam}{15}\newcommand{\sqrad}{7.5}
\newcommand{\bigdiam}{18}\newcommand{\bigrad}{9}
\newcommand{\width}{53}
\newcommand{\bigsqdiam}{18}\newcommand{\bigsqrad}{9}
\newcommand{\bigsqwidth}{35}
\newcommand{\psquare}[1]{\makebox(0,0){$#1$}%
\put(-\sqrad,-\sqrad){\line(1,0){\sqdiam}\line(0,1){\sqdiam}}\put(\sqrad,\sqrad){\line(-1,0){\sqdiam}\line(0,-1){\sqdiam}}}
\newcommand{\poval}[1]{\makebox(0,0){$#1$}\oval(20,20)}
\newcommand{\ppoval}[3]{%
\put(0,0){\makebox(0,0){$#1$}}\put(#3,0){\makebox(0,0){$#2$}}%
\oval(\bigdiam,\bigdiam)[l]\put(0,\bigrad){\line(1,0){#3}}\put(0,-\bigrad){\line(1,0){#3}}\put(#3,0){\oval(\bigdiam,\bigdiam)[r]}}
\begin{picture}(0,60)
%\put(-55,0){%
%\putmathlabel{-25}{50}{\psquare{A}}
%\drawline(-26,44)(-36,10)
%\drawline(-22,44)(12,10)
%\putmathlabel{-39}{4}{\ppoval{A}{A}{\width}}%
%}
\put(0,0){%
\putmathlabel{-25}{44}{\poval{A}}
\drawline(-24,44)(-24,10)
\drawline(24,44)(24,10)
\putmathlabel{-24}{4}{\ppoval{A}{C}{48}}}
\end{picture}
\end{equation}
%\end{center}
%
\begin{figure*}[t]
\newcommand{\boxwidth}{30}
\newcommand{\twoboxwidth}{60}
\newcommand{\boxradx}{15}
\newcommand{\boxrady}{6}
\newcommand{\boxheight}{12}
\newcommand{\layer}{42}
\newcommand{\drawbox}[1]{\put(0,0){\line(1,0){\boxwidth}\line(0,1){\boxheight}}%
\put(\boxwidth,\boxheight){\line(-1,0){\boxwidth}\line(0,-1){\boxheight}}\put(\boxradx,\boxrady){\makebox(0,0){$#1$}}}
\newcommand{\boxline}[1]{\put(\boxradx,\boxheight){\line(#1,1){\boxwidth}}}
\newcommand{\boxlinel}[3]{\boxline{-1}\put(0,\boxheight){\put(#2,#3){\putmathlabel{0}{\boxradx}{#1}}}}
\newcommand{\boxlinem}[3]{\boxline{0}\put(\boxradx,\boxheight){\put(#2,#3){\putmathlabel{0}{\boxradx}{#1}}}}
\newcommand{\boxliner}[3]{\boxline{1}\put(\boxwidth,\boxheight){\put(#2,#3){\putmathlabel{0}{\boxradx}{#1}}}}
\newcommand{\partlabel}[1]{\put(\boxwidth,-15){\put(\boxradx,0){\makebox(0,0){(#1)}}}}
\begin{center}
\begin{picture}(0,120)(0,-15)
\put(-220,0){%
 \partlabel{a}
 \put(\boxwidth,0){\drawbox{C}\put(\boxradx,0){\thicklines\boxlinel{g}{-5}{-5}}}
 \put(\twoboxwidth,0){\drawbox{\ucut}\boxlinem{\id}{7}{0}}
 \put(0,\layer){%
  \put(0,0){\drawbox{\uax}\boxlinem{\id}{-7}{0}}
  \put(\boxwidth,0){\drawbox{B}\put(\boxradx,0){\thicklines\boxlinel{f}{8}{4}}}
  \put(\twoboxwidth,0){\drawbox{\ucut}}
  \put(0,\layer){%
   \put(0,0){\drawbox{\uax}}
   \put(\boxwidth,0){\drawbox{A}}
  }
 }
}%
\put(-50,0){\put(-\boxradx,0){%
 \partlabel{b}
  \drawbox{C}\put(\boxradx,0){\thicklines\boxlinem{g}{-7}{0}}%
  \put(\boxwidth,0){\drawbox{\gamma}}
  \put(\twoboxwidth,0){\drawbox{\beta}\boxlinem{\id}{7}{0}}
  \put(0,\layer){%
   \drawbox{b}\boxlinem{\id}{-7}{0}%
   \put(\boxwidth,0){\drawbox{B}\put(\boxradx,0){\thicklines\boxlinem{f}{7}{0}}}
   \put(\twoboxwidth,0){\drawbox{\beta}}
   \put(0,\layer){%
    \drawbox{b}%
    \put(\boxwidth,0){\drawbox{a}}
    \put(\twoboxwidth,0){\drawbox{A}}
  }
 }
}}%
%
% OLD 5-LAYER VERSION INVOLVING TWIST
%\put(-\boxwidth,0){\partlabel{b}%
% \drawbox{C}\boxline{0}%
% \put(\boxwidth,0){\drawbox{\beta}\boxline{1}}
% \put(\twoboxwidth,0){\drawbox{\gamma}\boxline{-1}}
% \put(0,\layer){%
%  \drawbox{C}\put(\boxradx,0){\thicklines\boxlinem{g}{-7}{0}}%
%  \put(\boxwidth,0){\drawbox{\gamma}}
%  \put(\twoboxwidth,0){\drawbox{\beta}\boxline{0}}
%  \put(0,\layer){%
%   \drawbox{b}\boxline{0}%
%   \put(\boxwidth,0){\drawbox{B}\put(\boxradx,0){\thicklines\boxlinem{f}{7}{0}}}
%   \put(\twoboxwidth,0){\drawbox{\beta}}
%   \put(0,\layer){%
%    \drawbox{b}\boxline{1}%
%    \put(\boxwidth,0){\drawbox{a}\boxline{-1}}
%    \put(\twoboxwidth,0){\drawbox{A}\boxline{0}}
%    \put(0,\layer){%
%     \drawbox{a}%
%     \put(\boxwidth,0){\drawbox{b}}
%     \put(\twoboxwidth,0){\drawbox{A}}
%    } 
%  }
%  }
% }
%}
%
\newcommand{\loopp}[1]{\qbezier(0,0)(7,#1)(14,0)}
\newcommand{\looppp}[1]{\qbezier(0,0)(14,#1)(28,0)}
\put(95,0){%
 \putv{0}{0}\putv{14}{0}\putv{28}{0}\putv{42}{0}
 \put(14,0){\line(0,1){\layer}}\put(28,0){\line(0,1){\layer}}
 \put(28,0){\loopp{14}}
 \put(5,19){$g$}
 \put(0,\layer){%
  \putv{0}{0}\putv{14}{0}\putv{28}{0}\putv{42}{0}
  \put(0,0){\line(0,1){\layer}\loopp{-14}}\put(14,0){\loopp{14}}
  \put(28,0){\loopp{-14}}\put(42,0){\line(-1,3){14}}
  \put(14,18){$f$}
  \put(0,\layer){%
   \putv{0}{0}\putv{14}{0}\putv{28}{0}\putv{42}{0}
   \put(28,0){\loopp{-14}}
  } 
 }
 \put(17,0){\partlabel{c}}
 \put(80,0){
  \putv{0}{0}\putv{14}{0}\putv{28}{0}\putv{42}{0}
  \put(14,0){\line(1,6){14}}\put(28,0){\line(-1,3){28}}
  \put(28,0){\loopp{14}}
  \put(0,\layer){%
   \put(32,0){$f;g$}
   \put(0,\layer){%
    \putv{0}{0}\putv{14}{0}\putv{28}{0}\putv{42}{0}
    \put(0,0){\looppp{-20}}
   \put(28,0){\loopp{-14}}
   } 
  }
 }
}
\end{picture}
\end{center}
\caption{\label{three}Three approaches to incorporating axioms $p\vee\dual{p}$.  The diagrams above illustrate
composition schematically. (a) Define an abstract classical proof
$f:A\to B$ as a $\lgrel{L}$ morphism $f:\uax\wedge A\to
B\vee\ucut$, where the (potentially infinite) \emph{universal
axiom} $\uax$ is the product of all axioms $p\vee\dual{p}$ and the
\emph{universal cut} is its dual $\ucut=\neg{\uax}$, the sum of all
cuts $\dual{p}\wedge p$.  (b) Define an abstract classical proof $A\to
B$ as a $\lgrel{L}$ morphism $a\wedge A\to B\vee\beta$, where $a$ is
any product of axioms $p\vee \dual{p}$ and $\beta$ is any sum of cuts
$\dual{p}\wedge p$.  In (a) and (b), linear distributivity is hidden
at the interface layer. (c) Follow the usual linear logic pattern and
relax the $\Rel$ morphisms of $\lgrel{L}$ to $\Link$ morphisms, where
$\Link$ is the category of sets and linkings.
%, a linking $X\to Y$ being a simple graph on $X+Y$.  
Composition 
%in $\Link$ 
is the usual alternating path composition (\ie, the `smooth' paths in the example above).%
%one (alternating or `smooth' paths), 
%as in Kelly-Mac Lane graphs.
}\end{figure*}%
and uniqueness fails in the MALL resolution condition \resprime\ since
there are two edges between an output resolution and an input
coresolution:
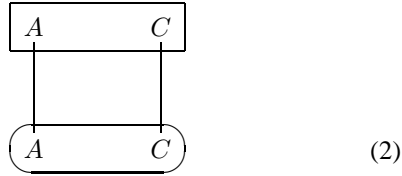
\begin{equation}\label{secondfailure}
\newcommand{\sqdiam}{15}\newcommand{\sqrad}{7.5}
\newcommand{\bigdiam}{18}\newcommand{\bigrad}{9}
\newcommand{\width}{53}
\newcommand{\bigsqdiam}{18}\newcommand{\bigsqrad}{9}
\newcommand{\bigsqwidth}{35}
\newcommand{\psquare}[1]{\makebox(0,0){$#1$}%
\put(-\sqrad,-\sqrad){\line(1,0){\sqdiam}\line(0,1){\sqdiam}}\put(\sqrad,\sqrad){\line(-1,0){\sqdiam}\line(0,-1){\sqdiam}}}
\newcommand{\poval}[1]{\makebox(0,0){$#1$}\oval(20,20)}
\newcommand{\ppoval}[3]{%
\put(0,0){\makebox(0,0){$#1$}}\put(#3,0){\makebox(0,0){$#2$}}%
\oval(\bigdiam,\bigdiam)[l]\put(0,\bigrad){\line(1,0){#3}}\put(0,-\bigrad){\line(1,0){#3}}\put(#3,0){\oval(\bigdiam,\bigdiam)[r]}}
\newcommand{\ppsquare}[3]{%
\put(0,0){\makebox(0,0){$#1$}}\put(#3,0){\makebox(0,0){$#2$}}%
\sqleft{\bigrad}\put(0,\bigrad){\line(1,0){#3}}\put(0,-\bigrad){\line(1,0){#3}}\put(#3,0){\sqright{\bigrad}}}
%\put(0,0){\makebox(0,0){$#1$}}\put(#3,0){\makebox(0,0){$#2$}}%
%\put(-\bigsqrad,-\bigsqrad){\line(1,0){#3}\line(0,1){\bigsqdiam}}%
%\put(-\bigsqrad,\bigsqrad){\line(0,-1){\bigsqdiam}\line(1,0){#3}}}
\begin{picture}(0,60)
%\put(-55,0){%
%\putmathlabel{-25}{50}{\psquare{A}}
%\drawline(-26,44)(-36,10)
%\drawline(-22,44)(12,10)
%\putmathlabel{-39}{4}{\ppoval{A}{A}{\width}}%
%}
\put(0,0){%
\putmathlabel{-24}{50}{\ppsquare{A}{C}{48}}
\drawline(-24,44)(-24,10)
\drawline(24,44)(24,10)
\putmathlabel{-24}{4}{\ppoval{A}{C}{48}}}
\end{picture}
\end{equation}
(Both failures depicted above apply to both distributions.)
These failures suggest naively relaxing the resolution conditions, to
admit distribution:
%Correspondingly, we relax the resolution conditions in the completely
%naive manner suggested by the distributions:
\begin{itemize}
\item[\lres] \emph{Lax resolution condition.} $R$ pulls resolutions back to super-resolutions and
pushes coresolutions forwards to super-coresolutions.
\end{itemize}
Here a super-(co)\-resolution is a superset of a (co)\-resolution.
Thus we have relaxed \res\ in the obvious way, admitting the first failure
(\ref{firstfail}) depicted above, by allowing the images of
(co)resolutions to spill beyond (co)resolutions.
Similarly, we relax \resprime\ in the obvious way to admit the
second failure (\ref{secondfailure}) depicted above, by simply dropping
uniqueness:\begin{itemize}
\item[\lresprime] \emph{Lax resolution condition.} $R$ has an edge between every 
output resolution and input coresolution.
\end{itemize}
Just as \resprime\ coincided with \res\ on syntactic
propositions, so \lresprime\ coincides with \lres.

We define the \emph{lax double glued category} $\glrel{K}$ using
%the lax
\lres\ in place of
%the original
\res.  
%It is not too hard to prove that, 
%
%We prove:
%
Surprisingly, this completely naive relaxation of the resolution
conditions, stimulated by the distribution failures, works:
\begin{itemize}
\item[{}] \emph{Theorem.} 
$\lgrel{K}$ contains as a full subcategory the free distributive
lattice category generated by the set $K$.
%The full subcategory $\dlk$ of $\G(\Rel^K)$ whose objects are
%syntactic propositions is the free distributive lattice category
%generated by the set $K$.
\end{itemize}
Do\u{s}en and Petri\'c \cite{DP04} define a distributive lattice category as a
binary product-sum category with a distribution, satisfying certain coherence
laws.  Thus 
%$\dlk$ is 
we obtain an abstract, syntax-free formulation of
contraction-weakening-distribution logic over a set of atoms $K$.

\paragraph*{Axioms.} The final step to an abstract, syntax-free formulation 
of classical propositional logic is to add axioms 
$1\to p\vee\overline{p}$ 
%to $\dl_L$ 
(hence by duality, also cuts $\overline{p}\wedge p\to 0$).
We explore three natural but distinct ways of achieving
this.
The first two constructions are quite general, and are not specific to
our $\Rel$-based abstract proposition approach, nor to the specific
axioms $p\vee \overline{p},q\vee\dual{q},\ldots$.  The third approach
is more ad hoc and limited, being $\Rel$-specific and
($p\vee\dual{p}$)-specific, but the style is more conventional in the
literature.
%
%In all three cases, the underlying preorder between abstract
%propositions is a boolean algebra.  
The three approaches are portrayed schematically in
Figure~\ref{three}.
\begin{itemize}
\myitem{a} \emph{Universal axiom construction.} Let the (potentially infinite)
abstract proposition $\uax$, the \emph{universal axiom}, be the
product of all $p\vee \overline{p}$ for complementary literals in $L$.
Its dual $\ucut=\neg \uax$ is the \emph{universal cut}. The
\emph{universal boolean category} $\bul$ has the objects of
$\lgrel{L}$ and a morphism $f:A\to B$ is a $\lgrel{L}$ morphism
\begin{center}\begin{picture}(0,40)\putmathlabel{0}{38}{\uax\wedge A}
\put(0,32){\vector(0,-1){21}}\putmathlabel{7}{20}{f}\putmathlabel{0}{3}{B\vee\ucut}\end{picture}\end{center}
Composition $f;g$ is defined in the obvious way, via linear
distribution $l$ at the interface.  See Figure~\ref{three}(a).
%\begin{center}\begin{picture}(0,120)
%\put(0,80){\putmathlabel{0}{38}{\uax\times A}\put(-4,21){\makebox(0,0)[r]{$\langle \pi_1, f\rangle$}}
%\put(0,32){\vector(0,-1){21}}\putmathlabel{0}{3}{\uax\times(B+\overline{\uax})}}
%\put(0,40){\put(0,34){\vector(0,-1){21}}\put(-4,23){\makebox(0,0)[r]{$l$}}\putmathlabel{0}{3}{(\uax\times B)+\overline{\uax}}}
%\put(0,0){\put(0,32){\vector(0,-1){21}}\put(-4,21){\makebox(0,0)[r]{$[ g,\inj_2]$}}\putmathlabel{0}{3}{C+\overline{\uax}}}
%\end{picture}\end{center}
\myitem{b} \emph{Local axiom construction.} The
\emph{local boolean category} $\bll$ has objects 
as above, but a morphism $A\to B$ is a $\lgrel{L}$ morphism $a\wedge
A\to B\vee\beta$ for $a$ a product of axioms and $\beta$ a sum of
cuts.  Composition is again defined in the obvious way, via linear
distribution.  See Figure~\ref{three}(b).
% and twist.
\myitem{c} \emph{Linkings.} We follow the standard recipe 
in linear logic and geometry of interaction \cite{Gir87}, imitating
the step from pure linearly distributive categories (two-sided proof
nets) to those with negation (one-sided nets) \cite{BCST96}.  The
\emph{linking boolean category} $\blinkingl$ is obtained from $\lgrel{L}$ 
by extending the homsets from $\Rel$ to the category $\Link$ of sets
and linkings.  A linking $X\to Y$ is a simple graph on $X+Y$, with
composition along alternating paths, like Kelly-Mac Lane graphs
\cite{KM71} (see Figure~\ref{three}(c)).
%  morphisms generalised from binary relations to
%linkings, which permit
%
%has 
%objects as above, and a morphism $A\to B$ is a binary relation on the
%leaves of
\end{itemize}
\paragraph*{Related work.} 
The category $\grel{K}$ was studied extensively in \cite{Hug04}, where
it was shown to fully embed the category of biextensional Chu spaces
over $K$ \cite{Bar79,Bar98}.  The definition of abstract proposition goes
via a tight orthogonality \cite{HS03} in $\grel{K}$, related to
totality spaces \cite{Loa94}.  The observation that the MALL
resolution condition \cite{HG03} characterises the free binary product-sum
category was in \cite{Hug02}\footnote{Read this technical report with
a pinch of salt, as the proof is far longer than it needs to be.  At
the time I wrote it, I was unaware of Hu's related work \cite{Hu99};
thus
\cite{Hug02} was never published.  Thanks to Robin Cockett and Robert Seely 
for pointing out the relationship with Hu's work.} in a more syntactic
guise (via the deductive system in \cite{CS01}); prior to that, Hongde
Hu had already characterised the free binary product-sum category in a
similar manner, using $P_4$-free graphs (contractible coherence
spaces) \cite{Hu99}.  In a syntactic setting, \cite{LS05} also
observes that
%Lamarche and Stra\ss{}burger work with the MALL
%resolution condition in \cite{LS05} in a syntactic setting, similarly
relaxing uniqueness in the MALL resolution condition yields a
classical proof net.  The classical proof nets sketched in
\cite{Gir91} are fleshed out in
\cite{Rob03}.  An abstract notion of classical proof 
net is presented in \cite{Hug04b}.  Categorical generalisations of
boolean algebras are presented in \cite{FP04} and \cite{DP04}.

\section{Abstract propositions}

Let $(X,S)$ be a set system, \ie, a set $X$ and a set $S$ of subsets
of $X$.  Subsets $s,t\subseteq X$ are \defn{orthogonal}, denoted
$s\,\bot\,t$, if
% $|s\cap t|=1$, \ie, 
they intersect in a single point.  The \defn{orthogonal of}
%$(X,S)\perp$ of $(X,S)$ is
%$(X,S\perp)$ where
$S$ is
\begin{displaymath}
S\perp \;\;=\;\; \{ t\subseteq X\;:\; t\,\bot\, s\text{ for all }s\in
S\}
\end{displaymath}
Fix a set of literals $L=\{p,\dual{p},q,\dual{q},\ldots\}$.
\begin{definition}
An \defn{abstract proposition} $(X,S)$ is a set $X$ of \defn{leaves},
each labelled by a literal, and a set $S$ of subsets of $X$, called
\defn{resolutions},
%pseudo-proposition $(X,S)$ 
%such that $(X,S)\perp{}\perp=(X,S)$.
satisfying:
\begin{itemize}
\item \emph{Double orthogonal:}  $S\perp{}\perp=S$.
%$(X,S)\perp{}\perp=(X,S)$.
%\item \emph{Independence.} For all $s,t\in S$, not $s\subseteq t$.
\end{itemize}
\end{definition}
%Fix a set of literals $L=\{p,\dual{p},q,\dual{q},\ldots\}$.  A
%\defn{pseudo-proposition} $(X,S)$ is a set $X$ of \defn{leaves},
%each labelled by a literal, and a set $S$ of subsets of $X$, called
%\defn{resolutions}.  Subsets $s,t\subseteq X$ are \defn{orthogonal},
%denoted $s\,\bot\,t$, if $|s\cap t|=1$, \ie, they intersect in a
%single leaf.  The \defn{orthogonal} 
%%$(X,S)\perp$ of $(X,S)$ is
%%$(X,S\perp)$ where
%of $S$ is
%\begin{displaymath}
%S\perp \;\;=\;\; \{ t\subseteq X\;:\; t\,\bot\, s\text{ for all }s\in
%S\}
%\end{displaymath}
%\begin{definition}
%An \defn{abstract proposition} is a pseudo-proposition $(X,S)$ 
%%such that $(X,S)\perp{}\perp=(X,S)$.
%satisfying:
%\begin{itemize}
%\item \emph{Double orthogonal.}  $S\perp{}\perp=S$.
%%$(X,S)\perp{}\perp=(X,S)$.
%%\item \emph{Independence.} For all $s,t\in S$, not $s\subseteq t$.
%\end{itemize}
%\end{definition}
%
Every syntactic $\wedge\vee$-formula $\phi$ over the set of literals
$L$ (\eg\ $(p\vee q)\wedge (p\vee \dual{p})$) defines an abstract
formula $(X,S)$ with $X$ the leaves of $\phi$: let $(X,E)$ be the
simple graph with an edge $xy\in E$ iff the leaves $x$ and $y$ meet at
a $\wedge$ in the parse tree of $\phi$, and let $S$ be the set of
maximal stable sets of $(X,E)$.  (A stable set of a graph is a maximal
set of vertices which contains no edge.) This represents $\phi$ modulo
associativity and commutativity of $\wedge$ and $\vee$. See page~1 for
an example. Any abstract proposition so obtained is \defn{syntactic}.

We define the following constants and operations.
%, in addition to $(-)\perp$ already defined.
\begin{itemize}
\item \emph{True}. $1=(\emptyset,\emptyset)$, no leaves and no resolutions.
\item \emph{False}. $0=(\emptyset,\{\emptyset\})$, no leaves and the empty resolution.
\item \emph{Negation/not}. $\neg(X,S)=(\overline{X},S\perp)$, where 
$\overline{X}$ relabels positive literals $p$ to negative literals
$\dual{p}$, and vice versa.
\item \emph{Sum/union/or}. %Given $(X,S)$ and $(Y,T)$, define 
$$\!\!(X,S)\vee (Y,T)\,=\,(X+Y,\,\{s+t:s\in S,t\in
T\}).$$ Thus the sum $A\vee B$ takes disjoint union on leaves, and a
resolution of $A\vee B$ is the union of a resolution of $A$ and a
resolution of $B$.
\item \emph{Product/join/and}.
 $A\wedge B=\neg(\neg A\vee \neg B)$.
% $A\wedge B=(A\perp\vee B\perp)\perp$.
\end{itemize}
%The operations $\wedge$ and $\vee$ coincide with on syntactic propositions, 
Note that $\neg\neg A=A$, $A\wedge 1=A$ and $A\vee 0=A$.

\paragraph*{Abstract truth.} An abstract
proposition is \defn{true} if every resolution contains a
complementary pair of leaves (\eg\ $p$ and $\dual{p}$).  Thus the
constant $1$ defined above (depicted on page~1)
% $\trueconst$
%depicted above 
is true (since it has no resolution) but
$0$
% $\falseconst$ 
is not (since it has an empty resolution).
%
%On non-empty abstract propositions this definition 
This definition of truth simply extends a well-known characterisation
of truth for syntactic formulas: a
%the resolutions of a 
syntactic formula $A$ is true iff every component of
%correspond to the components (disjuncts) of
its conjunctive normal form\footnote{The result of exhaustively
applying (co)distribution $A\vee(B\wedge C)\to (A\vee B)\wedge(A\vee C)$.}
$\textsc{cnf}(A)$ contains a complementary pair, and resolutions of
$A$ are in bijection with components of $\textsc{cnf}(A)$.
%, and $\textsc{cnf}(A)$
%is true iff each disjunct contains a complementary pair.
%%
%The resolutions of a syntactic proposition $A$ correspond to the
%components (disjuncts) of its conjunctive normal form
%$\textsc{cnf}(A)$.  For example, $A=p\vee (q\wedge \overline{p})$ has
%two resolutions,
%\begin{center}
%\begin{picture}(0,60)
%\thicklines\put(-35,30){\putmathlabel{-60}{0}{p\vee (q\wedge \overline{p})}
%\putmathlabel{-15}{0}{\mapsto}\put(30,0){\veegraph}\putmathlabel{72}{0}{\mapsto}\put(125,0){\veeresolutions}}
%\end{picture}
%\end{center}
%and $\textsc{cnf}(A)=(p\vee q)\wedge(p\vee\overline{p})$, whose two
%disjuncts $(p\vee q)$ and $(p\vee\overline{p})$ correspond to the two
%resolutions of $A$.  Thus this abstract definition of truth coincides
%with the usual syntactic one (since $A$ is true iff every disjunct of
%$\textsc{cnf}(A)$ contains a complementary pair of literals).
The abstract leaf/resolution representation of propositions can be
seen as ``\textsc{cnf} + superposition information'', where the latter
indicates which literal occurrences in different components of
$\textsc{cnf}(A)$ came from the same occurrence in $A$.

\section{Abstract proofs}

Abstract propositions correspond to objects of the category
$\G(\Rel^{L})$ where $\G$ is the Hyland-Tan double glueing
construction and $\Rel$ is the standard category of sets and
relations: $(X,S)$ corresponds to the $\G(\Rel^L)$ object
$(X,S\perp,S)$.  This category fully embeds the category of
biextensional Chu spaces over $L$ \cite{Hug04}.  Abstract propostions
correspond precisely to the objects of the tight orthogonality
subcategory, in the sense of
\cite{HS03}, for $s\,\bot\,t$ defined above.
The definitions of $0$, $1$, $\wedge$ and $\vee$ above correspond to
the standard initial object, terminal object, product and sum in a
tight orthogonality category.
%
%%Write $\trel{L}$ for the tight orthogonality subcategory of
%%$\grel{L}$, \ie, the category with abstract propositions as objects
%%and $\grel{L}$ morphisms between them.
%Write $\ap{L}$ for the category with abstract propostions as objects
%and $\grel{L}$ morphisms between them.  Thus $\ap{L}$ is isomorphic to
%the tight orthogonality subcategory of $\grel{L}$.

%The double orthogonality property implies that
%`cocoresolutions' coincide with resolutions.
%, \ie, a resolution of $\neg(X,S)$.  
The morphisms of $\grel{L}$
%(see \eg\ \cite{Hug04}) 
provide a ready-made notion of abstract proof.  
A \defn{coresolution} of an abstract proposition $(X,S)$ is an element
of $S\perp$.  
An \defn{abstract proof} $(X,S)\to (Y,T)$ between abstract
propositions is a binary relation $R$ between $X$ and $Y$ which
respects labelling, \ie, $xRy$ only if $x$ and $y$ are labelled with
the same literal, and satisfies:
\begin{itemize}
\item[\res] \emph{Resolution condition.} $R$ pulls resolutions backwards and
pushes coresolutions forwards.\footnote{For logical reaons, we have
taken resolutions as the contravariant part, so the ``co'' is opposite
to usual $\G(\mathcal{C})$. In the general $\G(\mathcal{C})$ case, the sets of
resolutions and coresolutions are independent, rather than the one
determining the other by orthogonality, as we have with abstract propositions.}
\end{itemize}
More precisely, the inverse image of $R$ is a function $T\to S$ and
the direct image of $R$ is a function $S\perp\to T\perp$.
As discussed in the Introduction, this coincides with the usual
resolution condition on MALL proof nets \cite{HG03}:
%\begin{itemize}
%\item \emph{Resolution condition.} $R$ has a unique edge on any resolution of $A\to B$.
%\end{itemize}
%Here a resolution of $A\to B$ is a resolution of $\overline{A}\vee B$,
%\ie, a union of a coresolution of $A$ and a resolution of $B$.
%In other words:
\begin{proposition}\label{coincide} Between syntactic propositions, the resolution condition \res\ on $R:A\to B$ 
coincides with:
\begin{itemize}
\item[\resprime] \emph{Resolution condition.} $R$ has a unique edge between any 
output resolution and input coresolution.
\end{itemize}
\end{proposition}
\begin{proof}
Suppose $R$ satisfies \res, let $b$ be a resolution of $B$, and let
$\alpha$ be a coresolution of $A$.  Let $a$ be the resolution of $A$ which
is the $R$-image of $b$.  By double orthogonality, $a$ and $\alpha$
intersect at a single leaf $x$.  Since $R$ maps $b$ onto $a$, we have
$xRy$ for some leaf $y$ in $b$.  This provides the unique edge between
$\alpha$ and $b$, for if $xRy'$ for some other $y\in b$, the
coresolution $R(\alpha)$ of $B$ would intersect $b$ in two leaves, a
contradiction.

Conversely, suppose $R$ satisfies \resprime, and let $b$ be a resolution
of $B$. Given a coresolution $\alpha$ of $A$ write $\widehat{\alpha}$
for the leaf of $A=(X,S)$ which is in the unique edge of $R$ between
$\alpha$ and $b$. Thus the $R$-image of $b$ is
$a=\{\widehat{\alpha}:\alpha\in S\perp\}$.  Since $a$ intersects each
$\alpha$ in exactly one leaf, namely $\widehat{\alpha}$, by the double
orthogonality condition $a$ is resolution.
\end{proof}
%An abstract proposition $(X,S)$ is \defn{finite} if $X$ is finite, and
%\defn{empty} if $S$ or $S\perp$ is empty.  
Let $\cwk$ be the full subcategory of $\grel{K}$ whose objects are
syntactic propositions.
\begin{theorem}\label{thm-cwk}
$\cwk$ is the free binary product-sum category generated by the set $K$.
\end{theorem}
\begin{proof}
By the Whitman-style theorem for binary product-sum categories in \cite{Hu99},
it essentially suffices to show that $\cwk$ is soft in the sense of
\cite{Joy95}: every morphism $A\wedge B\to C\vee D$ factors through a
projection on the left or an injection on the right.
If softness failed, there would be edges $A$--$C$ and $B$--$D$ (or
$A$--$D$ and $B$--$C$), breaking the resolution condition.
\end{proof}
Here we mean \emph{free category} in the standard 2-categorical sense
(see \eg\ \cite{Pow98}).  To obtain the stronger result, in which
$\cwk$ is initial, adjust the definition of \emph{syntactic} in the
obvious way: identify each element $k$ of $K$ with a vertex labelled
$k$, and define \emph{syntactic} to mean \emph{generated from such
vertices by $\wedge$ and $\vee$}.  (Thus the objects are in bijection
with $\wedge/\vee$-formulas over $K$.)

\section{Distribution forces lax resolution}

The Introduction discussed how the obvious candidates for distribution
fail the resolution condition.  We proceed completely naively, and
relax the resolution condition on a binary relation $R$ in the obvious
way to accomodate distribution.
\begin{itemize}
\item[\lres] \emph{Lax resolution condition.} $R$ pulls resolutions back to super-resolutions and
pushes coresolutions forwards to super-coresolutions.
\end{itemize}
Here a super-(co)\-resolution is a superset of a (co)\-resolution.
\begin{itemize}
\item[\lresprime] \emph{Lax resolution condition.} $R$ has an edge between every 
output resolution and input coresolution.
\end{itemize}
\begin{proposition}
Between syntactic propositions, the conditions \lres\ and \lresprime\
coincide.
\end{proposition}
\begin{proof}
Similar to the proof of Proposition~\ref{coincide}.
\end{proof}
Define the \defn{lax double glued category} $\glrel{K}$ using
%the lax
\lres\ in place of
%the original
\res.  (It is relatively easy to see that composition of 
binary relations preserves the lax resolution condition.)  This $\laxg$
is a general lax double glueing construction $\laxg(\mathcal{C})$ only
when the homsets of $\mathcal{C}$ are equipped with a suitable $\le$
relation.  In the case $\mathcal{C}=\Rel$ (or $\Rel^K$), $\le$ is the
inclusion order.  

The abstract propositions $0$ and $1$ remain initial and terminal in
$\lgrel{K}$, since $\lgrel{K}$ remains structured over $\Rel$.
\begin{proposition}
In the lax setting, $\wedge$ and $\vee$ continue to be product and sum
on abstract propositions, \ie, in the tight orthogonality subcategory
of $\lgrel{K}$.
\end{proposition}
\begin{proof}
Straightforward from the definition of the lax resolution condition.
\end{proof}

Let $\dlk$ be the full subcategory of $\lgrel{K}$ whose objects are
syntactic propositions.
A \emph{distributive lattice category} is a category with binary
products, sums and a distribution, satisfying certain coherence
conditions \cite{DP04}.
\begin{theorem}\label{dthm}
$\dlk$ is the free distributive lattice category generated by the set
$K$.
\end{theorem}
Here \emph{free category} is once again in the standard 2-categorical
sense.  (See the paragraph following the proof of Theorem~\ref{thm-cwk}.)
The proof uses two key factorisation lemmas.\vspace*{-1ex}
%Here is the key lemma for the proof of the theorem.
\begin{lemma}[Distribution factorisation]
\mbox{}\;\;Any $\dlk$ morphism $R:A\wedge(B\vee C)\to D$ factorises through
distribution $d$, \ie, there exists $R'$ such that $R$ equals
$$A\wedge(B\vee C)\stackrel{d}{\longto} (A\wedge B)\vee (A\wedge
C)\stackrel{R'}{\longto} D.$$
\end{lemma}
\begin{proof}
Coresolutions of $A\wedge (B\vee C)$ are in bijection with
coresolutions of $(A\wedge B)\vee (A\wedge C)$.  Thus the $R'$ induced
by $R$ in the obvious way (duplicating every edge from a leaf of $A$)
is well-defined, and factorises through distribution.
\end{proof}
Let \defn{mix} $m_{A,B}:A\wedge B\to A\vee B$ be the
composite $$A\wedge B\to A\wedge(Z\vee B)\to
(A\wedge Z)\vee (A\wedge B)\to A\vee B$$ of injection, distribution
and a two projections, for some $Z$.  Thus the binary relation of
$m_{A,B}$ is simply the identity between leaves.
\begin{lemma}[Mix-softness factorisation]
\mbox{}\;\;Any $\dlk$ morphism $R:A\to B$ with $A$ a pure product\footnote{\Ie,
$A$ is a product of one or more leaves, \ie, $A=(X,S)$ with $S$
comprising every singleton $\{x\}$ for $x\in X$.}
% of one or more leaves
%(\ie, $A=(X,S)$ and $S$ comprises every singleton $\{x\}$)
and $B$ a pure sum
% of one or more leaves 
is \defn{mix-soft}: unless $R$ is the identity on a single leaf, it
factorises through an injection, a projection, or mix.
\end{lemma}
\begin{proof}
If any leaf is not covered, $R$ factorises through an injection or a projection.  Otherwise, assuming $R$ is not the identity on a single
leaf, it factorises through mix at whichever of $A$ or $B$ contains
more than one leaf.
\end{proof}

\begin{proofof}{Theorem~\ref{dthm}}
Coherence of distributive lattice categories with respect to a
faithful functor to $\Rel$ was proved in \cite{DP04}, and $\Rel$ is
the underlying morphism category of $\lgrel{K}$.  Thus we need only
show that every morphism of $\dlk$ is canonical, \ie, generated from
the canonical maps defining a distributive lattice category.

Given $R:A\to B$, using the distribution factorisation lemma
%by factorising through distribution $U\wedge(V\vee
%W)\to (U\wedge V)\vee (U\wedge W)$ 
we may assume $A$ is in disjunctive normal form (a sum of products of
leaves), and since $\dlk$ has sums, we may further assume $A$ is a pure product of leaves.
Dually, we may assume $B$ is a sum of leaves.  Apply the mix-softness factorisation lemma.
%Next, since $\dlk$ has sums, we can further
%assume that $A$ has a single disjunct, \ie, a product of leaves.
%%$A=(X,S)$ with $S$ all singletons $\{x\}$.  
%
%If $B$ is a product, reduce it using the two
%projections; if $B=B_1\vee B_2$, and there is no edge to one of the
%$B_i$, factorise through the injection into the other $B_i$; if
%$B=B_1\vee B_2$, factorise through mix $m_{B_1,B_2}$.  Thus we can
%reduce $B$ to a single leaf, and we are left with a morphism $A\to B$
%from a product of leaves to a leaf.  The resolution condition implies there is a
\end{proofof}
Note that the above proof, in terms of the two factorisation lemmas,
amounts to a Whitman-style characterisation theorem for free
distributive lattice categories.\footnote{Hence the generality of
defining mix as a distribution composite, rather than directly as
the identity binary relation between leaves.} (See \cite{Hu99} for the
pure binary product-sum case.)

\section{Axioms}

So far, with the lax double glued category $\lgrel{L}$, we have an
abstract setting for contraction-weakening-distribution logic.  The
setting is canonical in the sense that it fully embeds the free distributive lattice category.
The category $\lgrel{L}$ is also equipped with a duality $\neg$, a
contravariant full and faithful functor (a de Morgan duality).
The final step to an abstract, syntax-free formulation of classical
propositional logic is to add axioms $1\to p\vee\overline{p}$
%to $\dl_L$ 
(hence by duality, also cuts $\overline{p}\wedge p\to 0$).

We explore three natural but distinct ways of achieving this.  Each
was discussed and motivated in the Introduction, and the idea behind
composition was sketched in Figure~\ref{three}.

The first two constructions are quite general, and are not specific to
our $\Rel$-based abstract proposition approach, nor to the specific
axioms $p\vee \overline{p},q\vee\dual{q},\ldots$.  The third approach
is more ad hoc and limited, being $\Rel$-specific and
($p\vee\dual{p}$)-specific, but the style is more conventional in the
literature, \eg\ \cite{KM71,Gir87,BCST96}.
%
%In all three cases, the underlying preorder between abstract
%propositions is a boolean algebra.

\subsection{Universal axiom construction}

Let the (potentially infinite) abstract proposition $\uax$, the
\defn{universal axiom}, be the product of all $p\vee \overline{p}$ for
complementary literals in $L$.  Its dual $\ucut=\neg \uax$ is the
\defn{universal cut}. The \defn{universal boolean category} $\bul$ has the 
objects of $\lgrel{L}$ and a morphism $f:A\to B$ is a $\lgrel{L}$
morphism\vspace{-2ex}
\begin{center}\begin{picture}(0,40)\putmathlabel{0}{38}{\uax\wedge A}
\put(0,32){\vector(0,-1){21}}\putmathlabel{7}{20}{f}\putmathlabel{0}{3}{B\vee\ucut}\end{picture}\end{center}
Composition $f;g$ is defined in the obvious way, via linear
distribution $l$ at the interface:
\begin{center}\begin{picture}(0,120)
\put(0,80){\putmathlabel{0}{38}{\uax\wedge A}\put(-4,21){\makebox(0,0)[r]{$\langle \pi_1, f\rangle$}}
\put(0,32){\vector(0,-1){21}}\putmathlabel{0}{3}{\uax\wedge(B\vee\ucut)}}
\put(0,40){\put(0,34){\vector(0,-1){21}}\put(-4,23){\makebox(0,0)[r]{$l$}}\putmathlabel{0}{3}{(\uax\wedge B)\vee\ucut}}
\put(0,0){\put(0,32){\vector(0,-1){21}}\put(-4,21){\makebox(0,0)[r]{$[ g,\inj_2]$}}\putmathlabel{0}{3}{C\vee\ucut}}
\end{picture}\end{center}
On syntactic propositions, there is a morphism $A\to B$ iff
$A\Rightarrow B\,=\,\neg A\vee B$ is true.  An \defn{abstract
classical proof} of $A$ in $\bul$ is a morphism $1\to A$.  Thus $A$
has an abstract classical proof in $\bul$ iff it is true.

%The universal axiom construction is quite general.  Writing the
%product pairing $\langle h,k\rangle:U\to V\wedge W$ of $h:U\to V$ and
%$k:U\to W$ and sum pairing $[h,k]:U\vee V\to W$ of $h:U\to W$ and
%$k:V\to W$ with explicit adjunction units,
%\begin{center}\begin{picture}(0,200)(0,-40)
%\put(0,80){\putmathlabel{0}{38}{(\uax\wedge A)\wedge(\uax\wedge A)}\put(-4,21){\makebox(0,0)[r]{$\pi_1\wedge f$}}
%\put(0,32){\vector(0,-1){21}}\putmathlabel{0}{3}{\uax\wedge(B\vee\ucut)}
% \put(0,40){\put(0,32){\vector(0,-1){21}}\putmathlabel{0}{38}{\uax\wedge A}\put(-4,21){\makebox(0,0)[r]{$\delta$}}}}
%\put(0,40){\put(0,34){\vector(0,-1){21}}\put(-4,23){\makebox(0,0)[r]{$l$}}\putmathlabel{0}{3}{(\uax\wedge B)\vee\ucut}}
%\put(0,0){\put(0,32){\vector(0,-1){21}}\put(-4,21){\makebox(0,0)[r]{$g\vee\inj_2$}}\putmathlabel{0}{3}{(C\vee\ucut)\vee(C\vee\ucut)}}
%\put(0,-40){\put(0,32){\vector(0,-1){21}}\put(-4,21){\makebox(0,0)[r]{$\epsilon$}}\putmathlabel{0}{3}{C\vee\ucut}}
%\end{picture}\end{center}
%we see that we can apply the universal axiom construction to a category equipped with
%\begin{itemize}
%\item a tensor, with a linear distribution over its dual;
%\item a choice of axioms to be tensored together in forming $\uax$ (possibly an infinite tensor);
%\item indexed families of morphisms for contraction (typed like $\epsilon$), 
%weakening (typed like inclusion $\inj_i$), copying (typed like $\delta$) and
%deletion (typed like projection $\pi_j$);
%\item sufficient coherence laws to ensure associativity of the above composition by diagram 
%pasting.
%\end{itemize}

\subsection{Local axiom construction}

This construction is similar to the one above.
The \defn{local boolean category} $\bll$ has objects as above, but a
morphism $A\to B$ is a $\lgrel{L}$ morphism $a\wedge A\to B\vee\beta$
for $a$ a product of (zero or more) axioms $p\wedge \dual{p}$ and
$\beta$ a sum of (zero or more) cuts $\dual{p}\wedge p$.  Composition
is again defined in the obvious way, via associativity and linear
distribution at the interface:
\begin{center}\begin{picture}(0,205)(0,-40)
\put(0,80){%
 \put(0,40){\put(0,32){\vector(0,-1){21}}\putmathlabel{0}{38}{(b\wedge a)\wedge A}\put(-4,21){\makebox(0,0)[r]{$\assoc$}}}
  \putmathlabel{0}{38}{b\wedge (a\wedge A)}\put(-4,21){\makebox(0,0)[r]{$\id\wedge f$}}
  \put(0,32){\vector(0,-1){21}}\putmathlabel{0}{3}{b\wedge (B\vee\beta)}
}
\put(0,40){\put(0,34){\vector(0,-1){21}}\put(-4,23){\makebox(0,0)[r]{$l$}}\putmathlabel{0}{3}{(b\wedge B)\vee\beta}}
\put(0,0){\put(0,32){\vector(0,-1){21}}\put(-4,21){\makebox(0,0)[r]{$g\vee \id$}}\putmathlabel{0}{3}{(C\vee\gamma)\vee\beta}}
\put(0,-40){\put(0,32){\vector(0,-1){21}}\put(-4,22){\makebox(0,0)[r]{$\assoc$}}\putmathlabel{0}{3}{C\vee(\gamma\vee\beta)}}
\end{picture}\end{center}
See Figure~\ref{three}(b) for a schematic supressing canonical maps.
%In a manner similar to the universal axiom construction, the local
%axiom construction generalises to a category equipped with
%\begin{itemize}
%\item a tensor, with a linear distribution over its dual;
%\item a choice of axioms to be tensored together;
%\item indexed families of morphisms for contraction (typed like $\epsilon$), 
%weakening (typed like inclusion $\inj_i$), copying (typed like $\delta$) and
%deletion (typed like projection $\pi_j$);
%\item sufficient coherence laws to ensure associativity of the above composition by diagram 
%pasting.
%\end{itemize}
%%\section{Conclusion}
%In contrast to the universal construction, the local construction does
%not require the existence of an object corresponding to a possibly
%infinite tensor of axioms $p\vee\dual{p}$.

\subsection{Linkings}

Although the following approach is more ad hoc, in that it does not
generalise as the two constructions above, it is a standard idea in
linear logic \cite{Gir87}, traceable as far back as Kelly-Mac Lane
graphs for closed categories \cite{KM71}.  
%We follow the standard recipe in linear logic and
%geometry of interaction \cite{Gir87}, imitating the step from pure
%linearly distributive categories (two-sided proof nets) to those with
%negation (one-sided nets)
%\cite{BCST96}.  

Let $\Link$ denote the category of sets and linkings, where a linking
$X\to Y$ is a simple graph on the disjoint union $X+Y$.  Composition
is the usual alternating path composition (see
Figure~\ref{three}(c)).
%
%Let $\Link^0_L$ be the homset extension of $\lgrel{L}$ obtained by
%allowing arbitrary linkings between leaves which respect labelling.
%(Naturally, edges inside $A$ or inside $B$ must be between dual leaves
%($p$---$\dual{p}$.)
%
Let the category $\Link_L$ be the homset extension of $\lgrel{L}$
obtained by permitting arbitrary linkings between leaves which respect
labelling (with edges of $A\to B$ within $A$ or $B$ going between dual
leaves $p$---$\dual{p}$).  Since $\Link$ is compact closed under
disjoint union, $\Link_L$ is star-autonomous under $\wedge$.  Let
$\blinkingl$, the \defn{linking boolean category}, be the restriction
of $\Link_L$ to syntactic propositions while retaining the lax MALL
resolution condition \lresprime\ from $\lgrel{K}$, \ie, there is an
edge in every resolution of $A\Rightarrow B=\neg A\vee B$.  Since the
ambient category $\Link_L$ is star-autonomous, and objects are
syntactic, by a routine structural induction \lresprime\ is preserved
by composition.
It is immediate from \lresprime\ that $A\Rightarrow B$ is true iff
there is a morphism $A\to B$.
The lax MALL resolution condition is also studied in \cite{LS05}, in a
more syntactic setting.

%\footnotesize
\small
\bibliographystyle{alpha}
\bibliography{../bib/main}

\end{document}